\documentclass[12pt]{article}
\usepackage{amscd,amsfonts,amssymb,amsmath,amsthm,latexsym,mathrsfs,graphicx}
\usepackage[margin=1in]{geometry}

\newcommand{\alg}{\mathrm{alg}}
\newcommand{\an}{\mathrm{an}}

\newcommand{\bfD}{\mathbf{D}}

\newcommand{\bfG}{\mathbf{G}}

\newcommand{\bfP}{\mathbf{P}}
\newcommand{\bfQ}{\mathbf{Q}}

\newcommand{\bfZ}{\mathbf{Z}}

\newcommand{\bra}[1]{{\langle #1 \rangle}}
\newcommand{\bs}{\backslash}

\newcommand{\calL}{\mathcal{L}}

\newcommand{\calO}{\mathcal{O}}
\newcommand{\calR}{\mathcal{R}}

\newcommand{\cn}{\colon}

\newcommand{\crys}{\mathrm{crys}}

\newcommand{\de}{\delta}

\newcommand{\dR}{\mathrm{dR}}

\DeclareMathOperator{\Ext}{Ext}

\newcommand{\fkm}{\mathfrak{m}}

\newcommand{\from}{\leftarrow}

\newcommand{\Ga}{\Gamma}
\newcommand{\ga}{\gamma}
\DeclareMathOperator{\Gal}{Gal}

\DeclareMathOperator{\Hom}{Hom}

\newcommand{\inv}{^{-1}}

\newcommand{\la}{\lambda}

\newcommand{\llim}{\varprojlim}

\newcommand{\Om}{\Omega}

\DeclareMathOperator{\ord}{ord}

\newcommand{\pd}{\partial}

\newcommand{\rig}{\mathrm{rig}}

\newcommand{\st}{\mathrm{st}}

\newcommand{\vphi}{\varphi}

\newcommand{\pab}{\text{$p$-ab}}
\newcommand{\Qp}{{\bfQ_p}}
\newcommand{\Qpa}{{\bfQ_p^\alg}}
\newcommand{\Qpt}{{\bfQ_p^\times}}
\newcommand{\Qpth}{{\bfQ_p^{\times,\wedge}}}
\DeclareMathOperator{\wtt}{wt}
\newcommand{\Zp}{{\bfZ_p}}

\begin{document}

\hfil {\bf The $\calL$-invariant, the dual $\calL$-invariant, and
  families} \hfil

\hfil Jonathan Pottharst, draft \today \hfil

\hfil \emph{Dedicated to Glenn Stevens on the occasion of his 60th
  birthday} \hfil

\vskip 6pt

\hfil {\bf Abstract} \hfil

\vskip 6pt

\hfil \parbox{5.5in}{Given a rank two trianguline family of
  $(\vphi,\Ga)$-modules having a noncrystalline semistable member, we
  compute the Fontaine--Mazur $\calL$-invariant of that member in
  terms of the logarithmic derivative, with respect to the Sen weight,
  of the value at $p$ of the trianguline parameter.  This generalizes
  prior work, in the case of Galois representations, due to
  Greenberg--Stevens and Colmez.} \hfil

\vskip 12pt

\noindent {\bf \S0 Introduction}

In the remarkable paper [GS], Greenberg and Stevens proved a formula,
conjectured by Mazur, Tate, and Teitelbaum in [MTT], for the
derivative at $s=1$ of the $p$-adic $L$-function of an elliptic curve
$E/\bfQ$ when $p$ is a prime of split multiplicative reduction.  The
novel quantity in this formula was the so-called $\calL$-invariant,
namely $\calL(E) = \log_p(q_E)/\ord_p(q_E)$ where $q_E \in p\Zp$
generates the kernel of the Tate uniformization $\bfG_{m,\Qp}^\an \to
E_\Qp^\an$.  The proof of the Greenberg--Stevens Theorem had two main
steps.  On the one hand, a two-variable $p$-adic $L$-function was
constructed, allowing the sought derivative to be computed in terms of
the derivative of the Hecke $U_p$-operator with respect to the weight.
On the other hand, a local Galois cohomology computation was used to
relate this derivative of $U_p$ to the $\calL$-invariant, as in the
formula [GS, (0.15)].

This short paper extends [GS, (0.15)], by extending the technique of
its proof, to $(\vphi,\Ga)$-modules over the Robba ring.  The main
result, which is proved in \S2.3, is as follows.

\vskip 6pt

\noindent {\bf Theorem.}  \emph{Let $X$ be an analytic space over a
  finite extension $E$ of $\Qp$, let $\de,\eta \cn \Qp^\times \to
  \calO_X^\times$ be continuous characters, and let the
  $(\vphi,\Ga)$-module $D$ over $\calR_X$ be an extension of
  $\calR_X(\de)$ by $\calR_X(\eta)$.  Assume that $P \in X$ is such
  that the specialization $D_0$ of $D$ at $P$ is, up to twist,
  noncrystalline semistable.  Then the differential form}
\[
d\log(\eta\de\inv(p)) - \calL(D_0) \cdot d\wtt(\eta\de\inv)
  \in \Om_{X/E}
\]
\emph{vanishes at $P$.}

\vskip 6pt

We hope the extension might motivate the generalization of other
techniques of [GS] to eigenvarieties whose triangulations possess a
two-step graded piece that is noncrystalline semistable up to twist,
and not necessarily \'etale.  Although we had originally hoped to
study the Hodge--Tate property in a uniform manner by formulating
``dual $\calL$-invariants'', we show at the end of \S2.3 that such a
technique necessarily cannot work.

There is a rich history of notions of $\calL$-invariant, especially
for modular forms of higher weight, and comparisons among them, for
which we refer to [C].  Especially, the formula [GS, (0.15)] was
previously generalized to the nonordinary setting by Colmez [C2,
  Th\'eor\`eme~0.5], using a study of the adjoint representation and
delicate computations in Fontaine's rings.  Readers who are
comfortable translating between the languages of crystalline periods
and triangulations will find that his result is none other than the
special case of our result where the rank two family of
$(\vphi,\Ga)$-modules arises from a family of Galois representations.
(For other generalizations see [B], which we learned of after the
writing of this paper.)

For many years, Glenn Stevens has provided me with inspiration and
support.  It is a pleasure to dedicate this article to him.

\vskip 12pt

\noindent {\bf \S1 $\calL$-invariants in arithmetic}

\vskip 12pt

\noindent {\bf \S1.1 Abstract $\calL$-invariants}

Let $E$ be a field, $X$ a two-dimensional $E$-vector space with dual
$X^* = \Hom_E(X,E)$, and $\la^*,\mu^* \in X^*$ a distinguished ordered
basis.  For $x^* = a\la^* + b\mu^* \in X^*$, we write $\bra{\la^*,x^*}
= a$ and $\bra{\mu^*,x^*} = b$.

The \emph{$\calL$-invariant} is the bijection
\[
\calL \cn \bfP(X) \stackrel\sim\to \bfP^1(E),
\qquad
Ex \mapsto (\la^*(x) : \mu^*(x)),
\]
and the \emph{dual $\calL$-invariant} is the bijection
\[
\calL^* \cn \bfP(X^*) \stackrel\sim\to \bfP^1(E),
\qquad
Ex^* \mapsto (-\bra{\mu^*,x^*} : \bra{\la^*,x^*}).
\]
We often confuse $\calL$ (resp.\ $\calL^*$) with its composition with
the projection $X \bs \{0\} \to \bfP(X)$ (resp.\ $X^* \bs \{0\} \to
\bfP(X^*)$).

Note that two lines $L \subset X$ and $L^* \subset X^*$ are orthogonal
if and only if $\calL(L) = \calL^*(L^*)$.

\vskip 12pt

\noindent {\bf \S1.2 Arithmetic setup}

We fix a prime $p$, choose an algebraic closure $\Qpa$ of $\Qp$, and
let $G_\Qp = \Gal(\Qpa/\Qp)$.  The base field (not to be confused with
the coefficient field) of all Galois representations
(resp.\ $(\vphi,\Ga)$-modules) will always be $\Qp$ in this paper, so
we omit it from the notations of continuous Galois cohomology
(resp.\ Herr cohomology).

Write $\Qpth = \llim_n \Qpt/(\Qpt)^{p^n}$ for the pro-$p$ completion
of $\Qpt$.  Kummer theory gives rise to an identification $\Qpth
\stackrel\sim\to H^1(\Zp(1))$, and using flat base change for
$\otimes_\Zp S$ we deduce from this an isomorphism $\Qpth \otimes_\Zp
S \stackrel\sim\to H^1(S(1))$, whose inverse we denote by $q$.  On the
other hand, the local Artin map induces an isomorphism $\Qpth
\stackrel\sim\to G_\Qp^\pab$ (the maximal pro-$p$-abelian quotient),
and by composition an identification $H^1(S) = \Hom_\Zp(G_\Qp^\pab,S)
= \Hom_\Zp(\Qpth,S)$, which we denote by $q^*$.  The pairing
\[
H^1(S(1)) \times H^1(S)
\cong
(\Qpth \otimes_\Zp S) \times \Hom_\Zp(\Qpth,S)
\to S,
\]
given $q \times q^*$ followed by evaluation of homomorphisms,
coincides (up to sign) with the Tate pairing, given by cup product and
the local invariant map.  It is therefore perfect when $S$ is a finite
field extension of $\Qp$.

We denote by $\log_p \cn \Qpt \to \Qp$ the natural extension of the
logarithm satisfying $\log_p(p)=0$, and by $\ord_p \cn \Qpt \to \bfZ$
the $p$-adic valuation; each of these induces by continuity a
homomorphism $\Qpth \to \Qp$, which we denote by the same name. Under
$q^*$, the $S$-module $H^1(S)$ is free of rank two on the ordered
basis $\log_p,\ord_p$.

The statement of the theorem assumes the reader is familiar with the
language of $(\vphi,\Ga)$-modules $D$ over Berger's Robba ring
$\calR_X$ with coefficients in an analytic space $X$ over $\Qp$.
However, in our work below we will only need the case where $X$ is the
spectrum of a finite $\Qp$-algebra $B$, in which case simply $\calR_B
= \calR_\Qp \otimes_\Qp B$.  We will need to refer to Fontaine's $2\pi
i$-element $t$, the objects of character type $\calR_B(\de)$ (where
$\de \cn \Qpt \to B^\times$ is a continuous character), and their Herr
cohomology; see [KPX] for details.  We denote by $x,|x| \cn \Qpt \to
\Qpt$ the identity map and $p$-adic absolute value, respectively, when
we wish to emphasize them as continuous characters of $\Qpt$.  Then,
for example, $\calR_\Qp(x) \cong t\calR_\Qp \subset \calR_\Qp$ and
$\calR_\Qp(x\cdot|x|) \cong \bfD_\rig(\Qp(1))$.

\vskip 12pt

\noindent {\bf \S1.3 Extensions of characters and their
  $\calL$-invariants}

In this subsection $E/\Qp$ is a finite field extension, and we apply
system of \S1.1 to the data $X = H^1(E(1))$, $X^* = H^1(E)$, $\la^* =
\log_p$, $\mu^* = \ord_p$.  Fix also continuous characters
$\de_0,\eta_0 \cn \Qpt \to E^\times$, and a $(\vphi,\Ga)$-module $D_0$
over $\calR_E$ sitting in a short exact sequence
\[
0 \to \calR_E(\eta_0) \to D_0 \to \calR_E(\de_0) \to 0.
\]
This sequence defines an extension class $[D_0]$ in
\[
\Ext^1(\calR_E(\de_0),\calR_E(\eta_0))
=
\Ext^1(\calR_E,\calR_E(\eta_0\de_0\inv))
=
H^1(\calR_E(\eta_0\de_0\inv)),
\]
and knowledge of the span $E[D_0] \subseteq
H^1(\calR_E(\eta_0\de_0\inv))$ is equivalent to knowledge of the
isomorphism class of the $(\vphi,\Ga)$-module $D_0$ over $\calR_E$.
We restate the well-known classification of $D_0$ via (dual)
$\calL$-invariants.

\vskip 6pt

\emph{First case:} $\eta_0$ is not of the form $(x\cdot|x|)x^k\de_0$
or $x^{-k}\de_0$ for any integer $k \geq 0$.  Then one has $\dim_E
H^1(\calR_E(\eta_0\de_0\inv)) = 1$, so the nonsplit $D_0$ are all
isomorphic.

\vskip 6pt

\emph{Second case:} $\eta_0 = (x\cdot|x|)x^k\de_0$ for some integer $k
\geq 0$.  Thus $\calR_E(\eta_0) \cong t^k\calR_E(\de_0)(1)$, and one
has
\[
H^1(\calR_E(\eta_0\de_0\inv))
\cong
H^1(t^k\calR_E(1))
\stackrel\sim\to
H^1(\calR_E(1)) = H^1(E(1))
=
\Qpth \otimes_\Zp E.
\]
We write $q_{D_0}$ for the image of $[D_0]$ under this identification.
If $q_{D_0} \neq 0$, then knowledge of $D_0$ up to isomorphism is
equivalent to knowledge of its \emph{(Fontaine--Mazur)
  $\calL$-invariant} $\calL(D_0) = \calL(q_{D_0}) = (\log_p q_{D_0} :
\ord_p q_{D_0})$.

Such $D_0$ is semistable up to twist, and is moreover crystalline up
to twist if and only if $\ord_p q_{D_0} = 0$, that is, either $q_{D_0}
= 0$ or both $q_{D_0} \neq 0$ and $\calL(D_0) = \infty$.  Conversely,
whenever $D_0$ has rank two and is up to twist noncrystalline
semistable, we are in the above situation for uniquely determined
$\delta_0$ and $k$, and $\ord_p q_{D_0} \neq 0$ so that $\calL(D_0)$
is defined and $\calL(D_0) \neq \infty$.  Such $D_0$ are not
isomorphic to those arising in any other of these three cases.  In the
noncrystalline semistable case, a computation in Fontaine's theory
shows that if $e \in \bfD_\st(D_0(\de_0\inv))$ is a
$\vphi$-eigenvector mapping to a basis of $\bfD_\crys(\calR_E)$, then
$\calL(D_0)$ is the slope of the Hodge filtration on
$\bfD_\dR(D_0(\de_0\inv)) = \bfD_\st(D_0(\de_0\inv))$ relative to the
basis $e, N(e)$.  (There exist plenty rank two crystalline $D$ that
are not extensions of some $\calR_E(\de_0)$ by $t^kR_E(\de_0)(1)$.)

\vskip 6pt

\emph{Third case:} $\eta_0 = x^{-k}\de_0$ for some integer $k \geq 0$.
Thus $\calR_E(\eta_0) \cong t^{-k}\calR_E(\de_0)$, and one has
\[
H^1(\calR_E(\eta_0\de_0\inv))
\cong
H^1(t^{-k}\calR_E)
\stackrel\sim\from
H^1(\calR_E) = H^1(E) = \Hom_\Zp(\Qpth,E).
\]
We write $q^*_{D_0}$ for the image of $[D_0]$ under this
identification. If $q^*_{D_0} \neq 0$, then knowledge of $D_0$ up to
isomorphism is equivalent to knowledge of its \emph{(Hodge--Tate) dual
  $\calL$-invariant} $\calL^*(D_0) = \calL^*(q^*_{D_0}) =
(-\bra{\ord_p,q^*_{D_0}} : \bra{\log_p,q^*_{D_0}})$.

For example, all rank one objects are of the form $\calR_E(\de_0)$ for
uniquely determined $\de_0$, so the case where $k=0$ is none other
than the situation where $D_0$ is an extension of a general rank one
object by itself.  This extension is nonsplit if and only if
$q^*_{D_0} \neq 0$, in which case $D_0$ is Hodge--Tate up to twist if
and only if $\calL^*(D_0) = \infty$.

\vskip 12pt

\noindent {\bf \S2 $\calL$-invariants of specializations of families}

Throughout this section, we fix a finite extension $E/\Qp$, and a
first-order deformation $(B,\fkm)$ of $E$, to be defined immediately
below.

\vskip 12pt

\noindent {\bf \S2.1 First order deformations}

By a \emph{first-order deformation of $E$}, we mean an Artinian local
$E$-algebra $(B,\fkm)$ with $B/\fkm=E$ and $\fkm^2=0$, and by a
morphism of these we mean a local $E$-algebra map.  One has $B = E
\oplus \fkm$ as $E$-vector spaces, and $B^\times = E^\times \times
(1+\fkm)$.  Note that the K\"ahler derivative $d \cn B = E \oplus \fkm
\to \Om_{B/E}$ is zero on the first factor and an isomorphism $\fkm
\cong \Om_{B/E}$ on the second, and employing $d$ value-by-value with
respect to $\Qpth$ we deduce an identification $d \cn H^1(E) \otimes_E
\fkm = \Hom_\Zp(\Qpth,\fkm) \cong \Hom_\Zp(\Qpth,\Om_{B/E})$.

We may view $\bra{\log_p,\cdot}$, $\bra{\ord_p,\cdot}$ and $q \in
\Qpth$ as maps $H^1(E) \to E$, and thus also as maps $H^1(E) \otimes_E
\fkm \to \fkm$.  For $c \in H^1(E) \otimes_E \fkm$ we have
$\bra{\log_p,c}, \bra{\ord_p,c}, c(q) \in \fkm$, and
\[
\log_p \otimes \bra{\log_p,c} + \ord_p \otimes \bra{\ord_p,c} = c.
\]
Evaluating the preceding equation at $q=p$, we find that
$\bra{\ord_p,c} = c(p)$.  On the other hand, evaluating at any
nonidentity $\ga_0 \in 1+2p\Zp$, we find that $\bra{\log_p,c} =
\frac{c(\ga_0)}{\log_p \ga_0}$.

Let $c \in H^1(E) \otimes_E \fkm$ be given.  Writing $\fkm^* =
\Hom_E(\fkm,E)$, it is easy to see that the following conditions on
$c$ are equivalent:
\begin{itemize}
\item There exists nonzero $q \in H^1(E(1))$ such that $c \in
  H^1(E)^{q=0} \otimes_E \fkm$.
\item There exists nonzero $q \in H^1(E(1))$ such that $(q \otimes
  1)(c) = 0$.
\item One of $\bra{\log_p,c}, \bra{\ord_p,c} \in \fkm$ is an
  $E$-multiple of the other.
\item There exist $h \in H^1(E)$ and $m \in \fkm$ such that $c = h
  \otimes m$.
\item The $E$-subspace $L_c = \{(1 \otimes v)(c) \mid v \in \fkm^*\}
  \subseteq H^1(E)$ satisfies $\dim_E L_c \leq 1$.
\end{itemize}
If these conditions are satisfied (for example, whenever $\dim_E
\Omega_{B/E} = 1$) and moreover $c \neq 0$, we call $c$ a \emph{pure
  tensor}.  In this case $q$, $h$, and $m$ are uniquely determined up
to nonzero $E^\times$-multiples, one has $H^1(E)^{q=0} = L_c = Eh$,
and the lines $Eq \subseteq H^1(E(1))$ and $Eh \subseteq H^1(E)$ are
an orthogonal pair.  The orthogonality shows that $\calL(q) =
\calL^*(h)$, and this common quantity is the constant of
proportionality (with slight abuse of notation) $(-\bra{\ord_p,c} :
\bra{\log_p,c}) \in \bfP^1(E)$; we denote it by $\calL^*(c)$.  By
definition, if $\calL^*(c) \neq \infty$ one has
\begin{equation}
\bra{\ord_p,c} = -\calL^*(c)\bra{\log_p,c}.
\end{equation}

\vskip 12pt

\noindent {\bf \S2.2 Characters valued in first-order deformations}

A continuous character $\de \cn \Qpt \to B^\times$ can be written
uniquely in the form $\de = \de_0 \cdot (1+\de_1)$ where $\de_0 \cn
\Qpt \to E^\times$ and $\de_1 \cn \Qpth \to \fkm$ are continuous
homomorphisms.  We use subscripts to denote the formation of these
components; one has $(\de\eta)_0 = \de_0\eta_0$ and $(\de\eta)_1 =
\de_1 + \eta_1$, and for a morphism $f \cn B \to B'$ one has $(f \circ
\de)_0 = f \circ \de_0$ and $(f \circ \de)_1 = f \circ \de_1$.  One
has the short exact sequence
\[
0 \to \calR_E(\de_0) \otimes_E \fkm \to \calR_B(\de) \to
\calR_E(\de_0) \to 0,
\]
and the corresponding extension class in
\[
\Ext^1(\calR_E(\de_0),\calR_E(\de_0) \otimes_E \fkm)
=
H^1(\calR_E \otimes_E \fkm)
=
H^1(E) \otimes_E \fkm
\stackrel{q^*}=
\Hom_\Zp(\Qpth,\fkm)
\]
is computed by $\de_1$, which we often view via $(q^*)\inv$ as an
element of $H^1(E) \otimes_E \fkm$.  The sequence is thus nonsplit if
and only if $\de_1 \neq 0$, if and only if $\calR_B(\de)$ is not
isomorphic to $\calR_E(\de_0) \otimes_E B$.

We reinterpret the quantities of \S2.1 in differential language when
$c = \de_1$ for a continuous character $\de \cn \Qpt \to B^\times$,
assuming that $\de_1$ is a pure tensor.  We compute the image of
$\de_1$ under the identification $d \cn H^1(E) \otimes_E \fkm \cong
\Hom_\Zp(\Qpth,\Om_{B/E})$, value-by-value with respect to $\Qpth$, to
be
\[
d(\de_1)
= \frac{\de_0d(\de_1)}{\de_0}
= \frac{d(\de_0+\de_0\de_1)}{\de_0+\de_0\de_1}
= \frac{d(\de_0(1+\de_1))}{\de_0(1+\de_1)}
= d\log(\de),
\]
where the second step is because $\de_0$ (resp.\ $\de_1$) takes values
in $E$ (resp.\ $\fkm$).  In particular, evaluating both sides at $p$,
we find that
\[
d\bra{\ord_p,\de_1} = d(\de_1(p)) = d\log(\de(p)).
\]
Next we compute the image of
$d\bra{\log_p,\de_1} \in \fkm$ under the identification to be
\[
d\bra{\log_p,\de_1}
=
\frac{d\de_1(\ga_0)}{\log_p \ga_0}
=
\frac{d\log(\de(\ga_0))}{\log_p \ga_0}
=
-d\wtt(\de),
\]
where $\wtt(\de) \in B$ is the \emph{weight} of $\de$, defined as the
value of $-\frac{\log\de(\ga_0)}{\log_p\ga_0}$ for $\ga_0 \in 1+2p\Zp$
nonidentity and sufficiently close to the identity.  This invariant is
normalized so that $\wtt(\de)$ agrees with the Sen weight of
$\calR_B(\de)$, where we consider the Sen weight of $\Qp(1)$ to be
$-1$.  Substituting these calculations of $d\bra{\ord_p,\de_1}$ and
$d\bra{\log_p,\de_1}$ into the equation (1) above, we find that if
($\de_1$ is a pure tensor and) $\calL^*(\de_1) \neq \infty$ then
\begin{equation}
d\log(\de(p)) = \calL^*(\de_1) \cdot d\wtt(\de)
\text{ in } \Om_{B/E}.
\end{equation}

\vskip 12pt

\noindent {\bf \S2.3 Extensions of characters over first-order
  deformations}

We suppose given continuous characters $\de,\eta \cn \Qpt \to
B^\times$ and a short exact sequence
\[
0 \to \calR_B(\eta) \to D \to \calR_B(\de) \to 0,
\]
and set $D_0 = D \otimes_B E$.  We take Herr cohomology of the short
exact sequence
\[
0 \to \calR_E((\eta\de\inv)_0) \otimes_E \fkm
\to \calR_B(\eta\de\inv) \to
\calR_E((\eta\de\inv)_0) \to 0
\]
to obtain the exact sequence
\begin{equation}
H^1(\calR_B(\eta\de\inv))
\xrightarrow\alpha
H^1(\calR_E((\eta\de\inv)_0))
\xrightarrow\pd
H^2(\calR_E((\eta\de\inv)_0)) \otimes_E \fkm,
\end{equation}
in which the reduction map $\alpha$ sends the class of $D$ to $[D_0]$,
and the connecting map $\pd$ is given (up to sign) by cup product with
$(\eta\de\inv)_1 \in H^1(E) \otimes_E \fkm$.  In particular,
$(\eta\de\inv)_1 \cup [D_0] = 0$.

Now we treat separately the three cases of \S1.3 applied to $D_0$.

\vskip 6pt

\emph{Second case:} The Tate duality isomorphism
$H^2(\calR_E((\eta\de\inv)_0)) \otimes_E \fkm \cong \fkm$ corresponds
$(\eta\de\inv)_1 \cup [D_0]$ with $(\eta\de\inv)_1(q_{D_0})$ (up to
sign), so that $(\eta\de\inv)_1(q_{D_0}) = 0$.

\vskip 6pt

\noindent \emph{Proof of the theorem.}  By replacing $X$ by its
first-order infinitesimal neighborhood at $P$ and $E$ by the residue
field at $P$, we may assume $X$ is the spectrum of a first-order
deformation of $E$.  Recall that $D_0$ is noncrystalline semistable up
to twist, so that according to \S1.3 we are in the situation
immediately preceding this proof.  In particular, one has
$\ord_p(q_{D_0}) \neq 0$ and so $q_{D_0} \neq 0$, and also that
$\calL(D_0)$ is defined and $\calL(D_0) \neq \infty$.

If $(\eta\de\inv)_1 = 0$ then $d\log(\eta\de\inv(p)) =
d\wtt(\eta\de\inv) = 0$, and we are done.

If $(\eta\de\inv)_1 \neq 0$ then the identity
$(\eta\de\inv)_1(q_{D_0}) = 0$ implies that $(\eta\de\inv)_1$ is a
pure tensor, and by orthogonality $\calL(D_0) =
\calL^*((\eta\de\inv)_1)$.  The theorem now follows from (2) applied
to the character $\eta\de\inv$.  $\square$

\vskip 6pt

The two main ingredients in the proof of the theorem, namely the
computation (2) and the equation $(\eta\de\inv)_1(q_{D_0})=0$ coming
from the long exact sequence on local Galois cohomology, are also key
steps in the proof of [GS, (0.15)].

\vskip 6pt

\emph{First and third cases:} Now one has
$H^2(\calR_E((\eta\de\inv)_0)) = 0$, so the exactness of
(3) implies the reduction map $\alpha$ is surjective.

In the first case, we see that the unique split and nonspilt
possibilities for $D_0$ may occur.

In the third case, $[D_0]$ is identified to $q^*_{D_0} \in H^1(E)$.
The surjectivity of the reduction map $\alpha$ shows that the two
subspaces $Eq^*_{D_0} \subseteq H^1(E)$ and $E(\eta\de\inv)_1
\subseteq H^1(E) \otimes_E \fkm$ can be arbitrary of dimension at most
one.  Assuming that $[D_0] \neq 0$, so the dual $\calL$-invariant
$\calL^*(D_0)$ is defined, this independence of $q^*_{D_0}$ and
$(\eta\de\inv)_1$ shows that no formula for $\calL^*(D_0)$ purely in
terms of $\de,\eta$ is possible without further hypotheses.

\vskip 12pt

\noindent {\bf References}

[B] Denis Benois, Infinitesimal deformations and the
$\ell$-invariant.  \emph{Documenta Math.} Extra Volume: Andrei A.\
Suslin's Sixtieth Birthday (2010), 5--31.

[C] Pierre Colmez, Z\'eros suppl\'ementaires de fonctions $L$
$p$-adiques de formes modularies.  \emph{Algebra and Number Theory},
R.~Tandon (ed.), Hindustan Book Agency (2005), 193--210.

[C2] Pierre Colmez, Invariants $\calL$ et d\'eriv\'ees de valeurs
propres de Frobenius.  \emph{Ast\'erisque} {\bf 331} (2010), 13--28.

[GS] Ralph Greenberg and Glenn Stevens, $p$-adic $L$-functions and
$p$-adic periods of modular forms.  \emph{Invent.\ Math.}\ {\bf 111}
(1993), no.\ 2, 407--447.

[KPX] Kiran S.\ Kedlaya, Jonathan Pottharst, and Liang Xiao,
Cohomology of arithmetic families of $(\vphi,\Ga)$-modules.
\emph{J.\ Amer.\ Math.\ Soc.}\ {\bf 27} (2014), no.\ 4, 1043--1115.

[MTT] Barry Mazur, John Tate, and Jeremy Teitelbaum, On $p$-adic
analogues of the conjectures of Birch and Swinnerton-Dyer.
\emph{Invent\ Math.}\ {\bf 84} (1986), no.\ 1, 1--48.

[S] Glenn Stevens, Coleman's $\calL$-invariant and families of modular
forms.  \emph{Ast\'erisque} {\bf 331} (2010), 1--12.

\end{document}